\documentclass[conference]{IEEEtran}

\usepackage{amsmath,amssymb,amsfonts}
\usepackage{cite}

\ifCLASSINFOpdf
   \usepackage[pdftex]{graphicx}
\else
   \usepackage[dvips]{graphicx}
\fi

\usepackage{hyperref}
\hypersetup{
    colorlinks=false,
    linkcolor=blue,
    filecolor=magenta,      
    urlcolor=cyan,
    pdfpagemode=FullScreen,
    }

\usepackage{accents}
\usepackage{mdframed}
\usepackage{xcolor}

\newcommand{\R}{\mathbb{R}}

\newcommand{\vct}[1]{\boldsymbol{#1}}
\newcommand{\mtx}[1]{\boldsymbol{#1}}

\newcommand{\norm}[1]{\left|\left|#1\right|\right|}

\newcommand{\abs}[1]{\left|#1\right|}

\newcommand{\set}[1]{\mathcal{#1}}

\DeclareMathOperator*{\minimize}{\text{minimize}}
\DeclareMathOperator*{\maximize}{\text{maximize}}

\newcommand{\va}{\vct{a}}

\newcommand{\vp}{\vct{p}}

\newcommand{\vv}{\vct{v}}

\newcommand{\vx}{\vct{x}}

\newcommand{\mP}{\mtx{P}}

\newcommand{\setB}{\set{B}}

\newcommand{\setE}{\set{E}}

\newcommand{\setG}{\set{G}}

\newcommand{\setK}{\set{K}}

\newcommand{\setN}{\set{N}}

\newcommand{\setP}{\set{P}}

\newcommand{\setT}{\set{T}}

\usepackage{stmaryrd}
\usepackage{nomencl}
\usepackage{dsfont}

\makenomenclature

\renewcommand{\baselinestretch}{0.912}
\renewcommand{\IEEEbibitemsep}{0pt plus 0.5pt}
\makeatletter
\IEEEtriggercmd{\reset@font\normalfont\fontsize{7.9pt}{8.40pt}\selectfont}
\makeatother
\IEEEtriggeratref{1}

\begin{document}
\title{Managing Vehicle Charging During Emergencies via Conservative Distribution System Modeling}

\author{
\IEEEauthorblockN{Alejandro D. Owen Aquino, Samuel Talkington, Daniel K. Molzahn}
\IEEEauthorblockA{School of Electrical and Computer Engineering\\
Georgia Institute of Technology, Atlanta, GA\\
\{aaquino30, talkington, molzahn\}@gatech.edu}
}

\maketitle

\begin{abstract}
Combinatorial distribution system optimization problems, such as scheduling electric vehicle (EV) charging during evacuations, present significant computational challenges. These challenges stem from the large numbers of constraints, continuous variables, and discrete variables, coupled with the  unbalanced nature of distribution systems. In response to the escalating frequency of extreme events impacting electric power systems, this paper introduces a method that integrates sample-based conservative linear power flow approximations (CLAs) into an optimization framework. In particular, this integration aims to ameliorate the aforementioned challenges of distribution system optimization in the context of efficiently minimizing the charging time required for EVs in urban evacuation scenarios.

\end{abstract}


\nomenclature[01]{$\mathcal{E}$}{Set of lines}
\nomenclature[02]{$\mathcal{B}$}{Set of $j=1,\dots, \lvert \setB \rvert$ buses}
\nomenclature[03]{$\mathcal{P}_j$}{Set of phases $\rho \in \setP_j \subset \{a,b,c\}$ for bus $j$}
\nomenclature[04]{$\mathcal{N}$}{Set of $i=(j,\rho) \in  \setB \times \mathcal{P}_j$ single phase nodes}
\nomenclature[05]{$\mathcal{T}$}{Set of $t=1,\dots, T$ time steps for EV charging}
\nomenclature[06]{$\Xi$}{Set of transportation analysis zones (TAZs)}
\nomenclature[07]{$E_\xi$}{Set of $h=1,\dots, \lvert E_\xi \rvert$ EVs in TAZ $\xi$}
\nomenclature[08]{$\mathcal{K}$}{Set of buses $k \in \mathcal{K} \subset \setB$ with EVs in them}
\nomenclature[12]{$p_{EV_k}^t$}{Active EV charging demand at node $k$ at time $t$}
\nomenclature[15]{$v_i^t$}{Squared voltage magnitude at node $i$ at time $t$}
\nomenclature[16]{$\Gamma_{\sf max}$}{First time step when EVs start charging}
\nomenclature[17]{$\Lambda_{\sf max}$}{Maximum allowable sum for all voltage violations}
\nomenclature[18]{$\tau^t$}{Binary indicator if EV charging has begun by time $t$}
\nomenclature[19]{$C_\xi^t$}{Binary charging status of TAZ $\xi$ at time $t$}
\nomenclature[20]{$C_{\xi,h}^t$}{Binary charging status of EV $h$ in TAZ $\xi$ at time $t$}
\nomenclature[21]{$L_{\xi,h}^t$}{Battery level of EV $h$ in TAZ $\xi$ at time $t$}
\nomenclature[22]{$L_{\xi,h}^0$}{Starting battery level of EV $h$ in TAZ $\xi$}
\nomenclature[23]{$R$}{Charging rate of an EV}
\nomenclature[24]{$\beta$}{Time steps to fully charge an EV from 0\% at rate $R$}
\nomenclature[25]{$d_\xi$}{Departure time for TAZ $\xi$}
\nomenclature[26]{$x$,$\vx$,$\boldsymbol{X}$}{Scalar, vector, and matrix, respectively}
\nomenclature[27]{$\underline{x}$,$\overline{x}$}{Lower and upper bounds of $x$, respectively}
\nomenclature[28]{$\lambda^-$,$\lambda^+$}{Lower and upper bound slack variables, respectively}
\nomenclature[29]{$\rvert \cdot \lvert$}{Magnitude of a complex number or size of a set}
\nomenclature[30]{$\widetilde{x}$,$\undertilde{x}$}{Overestimate and underestimate of $x$, respectively}
\nomenclature[31]{$\mathds{1}$}{Vector of all ones}
\nomenclature[32]{$(\cdot)^\top$}{Transpose of a vector}


\section{Introduction}
Extreme weather events, wildfires \cite{abelmalak_enhancing_resilience_wildfires_2022,nazaripouya_grid_wildfire_2020}, and other natural disasters \cite{feng_can_we_evacuate_2020} have introduced significant challenges in power system operation. Broad efforts are being made to improve the resilience of power systems to these challenges \cite{panteli_influence_2015,abelmalak_enhancing_resilience_wildfires_2022}; however, they continue to pose  unsolved problems.

Particularly, the increasing penetration of electric vehicles (EVs) necessitates advanced evacuation scheduling techniques \cite{donaldson2022integration,donaldson_wildfire_ev_resilience_stochastic,abelmalak_enhancing_resilience_wildfires_2022}. Emergencies may come with inherently restrictive time horizons; in contrast, charging times for large fleets of EVs may be more flexible. Moreover, the resulting optimization requires binary variables to ensure consistent charging instructions to evacuating regions. These facts, combined with the unbalanced nature of distribution networks, cause it to be computationally intractable to incorporate the more realistic AC power flow equations into evacuation optimization. 

Therefore, we aim to solve these two problems jointly by leveraging a \emph{surrogate} linear approximation of a distribution network model estimated from samples of circuit quantities. In particular, we use an approximation that is \emph{conservative} \cite{pap_2021}, in the sense that the OPF constraint set is ensured to be satisfied by the approximation for the sampled points. Thus far, this technique has only been applied to transmission networks. This work is the first to extend it to distribution networks.

\subsection{Related work}
\subsubsection{Linear approximations of the power flow equations}
\label{sec:linear-approx}
Analytical and data-driven linearizations of the non-linear AC power flow equations are well-studied, with surveys in \cite{molzahn_hiskens-fnt2019} and \cite{jia_data_driven_approx_survey_2023}. Recognizing the limitations of conventional linearizations in security-critical contexts, we adopt the conservative linearization framework from \cite{pap_2021} that aims to minimize constraint violations in optimal power flow (OPF) problems.

\subsubsection{Evacuation of electric vehicles}
\label{sec:evacuation}
The escalation of extreme weather events motivates research into evacuation strategies and their impacts on EV charging \cite{feng_can_we_evacuate_2020, donaldson_integration_ev_resilience_2022, macdonald_ev_charging_capacity_evacuation_2021, hasan_large-scale_part1_2021, hasan_large-scale_part2_2021}. Existing research, employing stochastic models \cite{donaldson_wildfire_ev_resilience_stochastic} and sequential optimization \cite{donaldson2022integration}, has been limited to balanced transmission networks and small scales.

In contrast, our approach consolidates the optimization across the entire time horizon into a single problem, coordinating evacuation scheduling deterministically using binary variables, realistic zoning analyses as in \cite{hasan_large-scale_part1_2021, hasan_large-scale_part2_2021} and conservative sensitivity coefficient matrices as in \cite{pap_2021}. This method enables application to unbalanced distribution networks with a significantly expanded node count, enhancing scalability and applicability in comprehensive emergency scenarios.

\subsection{Contributions and paper outline}
\label{sec:intro:article-outline}

To the best of our knowledge, past power system literature has yet to consider the emergency EV charging problem at a level of fidelity consistent with the state of the art of general evacuation planning. In contrast with prior power systems EV charging literature, the general evacuation planning literature considers Transportation Analysis Zones (TAZs) \cite{hasan_large-scale_part1_2021,hasan_large-scale_part2_2021}. These TAZs are delegated at the level of towns and neighborhoods\textemdash thus, in practice, emergency EV charging problems should consider impacts on \emph{distribution systems}. This fact has been neglected by past literature on emergency EV charging, which only considers transmission-level analysis \cite{donaldson_integration_ev_resilience_2022,donaldson2022integration,macdonald_ev_charging_capacity_evacuation_2021}. Our prior work in~\cite{nilsson_owen_aquino_coogan_molzahn-GreenEVT} shows that evacuation schedules that do not consider electric vehicle charging can lead to substantial overloads of distribution system components.

In summary, the contributions of this research are threefold:
\begin{enumerate}
    \item An efficient and more realistic method to develop emergency charging strategies for electric vehicles at the level of TAZs\textemdash inherently a distribution systems problem.
    \item A method to embed estimated conservative linear distribution network models in optimization problems. 
    \item Increasing the computational tractability of the first contribution by synthesizing it with the second contribution.
\end{enumerate}

Section \ref{sec:citywide-evacuation-problem} presents an urban evacuation problem. Section \ref{sec:CLA-EEV-C} integrates conservative power flow linearizations and constraint generation into the problem. Section \ref{sec:numerical-experiments} provides numerical validation on a segment of the realistic distribution network of Greensboro, NC using the testbed in~\cite{nilsson_owen_aquino_coogan_molzahn-GreenEVT}.

\section{Urban Evacuation Problem}
\label{sec:citywide-evacuation-problem}
The urban evacuation scheduling problem aims to mobilize all residents to safety in the shortest possible time before the onset of a predictable natural disaster, such as a major hurricane \cite{feng_can_we_evacuate_2020} or wildfire \cite{donaldson_integration_ev_resilience_2022}. While scenarios involving solely internal combustion engine vehicles need only focus on modeling the city's transportation network, evacuation planners in cities with a high prevalence of EVs must account for the coupling between the transportation and electric distribution networks due to EV charging. Planners are then tasked with coordinating not only the departure times and routing of all neighborhoods, but also the charging instructions for all EVs to minimize overloading of the power system. 
%
%
%
%
\begin{figure}[t]
    \centering
    \includegraphics[width=\linewidth]{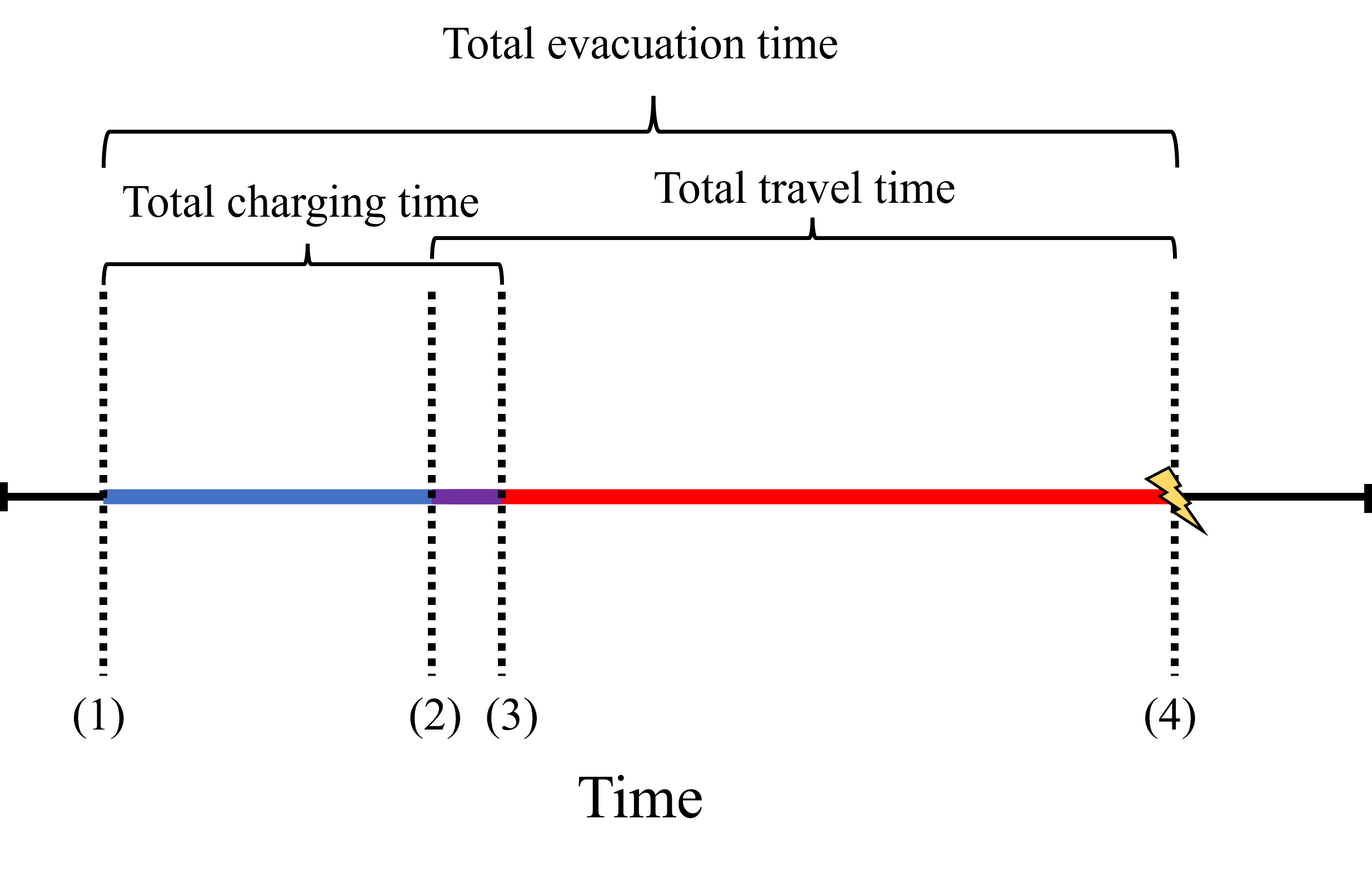}
    \caption{Illustration of an evacuation plan where charging and travel times are considered. The lightning symbol marks the beginning of a natural disaster.} 
    \label{fig: Evacuation_Horizon}
    \vspace*{-1.2em}
\end{figure}
A comprehensive evacuation timeline,  which illustrates both the charging and travel time components of the evacuation, is shown in Fig.~\ref{fig: Evacuation_Horizon}. This timeline is punctuated by four critical moments: (1) the start of EV charging, (2) the beginning of evacuation, (3) the end of EV charging, and (4) the completion of evacuation to safe locations prior to the disaster's onset.

Minimizing the entire timeline in Fig.~\ref{fig: Evacuation_Horizon} is a formidable challenge. As a result, this study focuses on optimizing the interval between moments (1) and (3) in the most effective manner. Developing an efficient formulation for the charging problem is a crucial step towards achieving our final objective: a systematic iterative process that combines solutions for the charging problem with the results for the departure-schedule-and-routing problem, as outlined in \cite{Huertas_2022}, to compute an optimal and comprehensive evacuation plan. 

\section{Challenges of the Emergency Electric Vehicle Charging (EEV-C) Problem}
\label{sec:challenges_and_methods}
We formulate and solve an emergency EV charging optimization problem in the framework of a disaster evacuation plan. Hereafter, we refer to this problem as the emergency EV charging problem (EEV-C). As in most distribution system optimization problems, the EEV-C problem is challenged by both the problem size and the presence of non-linear, non-convex engineering constraints. 

Additionally, solving the EEV-C problem presents a dilemma between (a) a longer charging schedule that avoids network violations and electrical component overloads, or (b) allowing violations to expedite the charging process. In the following subsections, we describe these challenges in more detail and propose methods to circumvent them. 

\subsection{Distribution network modeling}
\label{sec:distribution_opt_models}

Consider a distribution network model $\setG = (\setB,\setE)$, where the existence of a reference bus is implied. Each bus $j \in \setB$ can have up to three phases $\rho \in \setP_j \subset \{a,b,c\}$. Each single-phase \emph{node} is denoted as the tuple $i = (j,\rho) \in  \setN$. The size of $\setG$ adds substantial complexity on top of the fact that the power flow equations that govern $\setG$ are non-linear and non-convex. 

\begin{figure}[t]
    \centering
    \includegraphics[width=\linewidth]{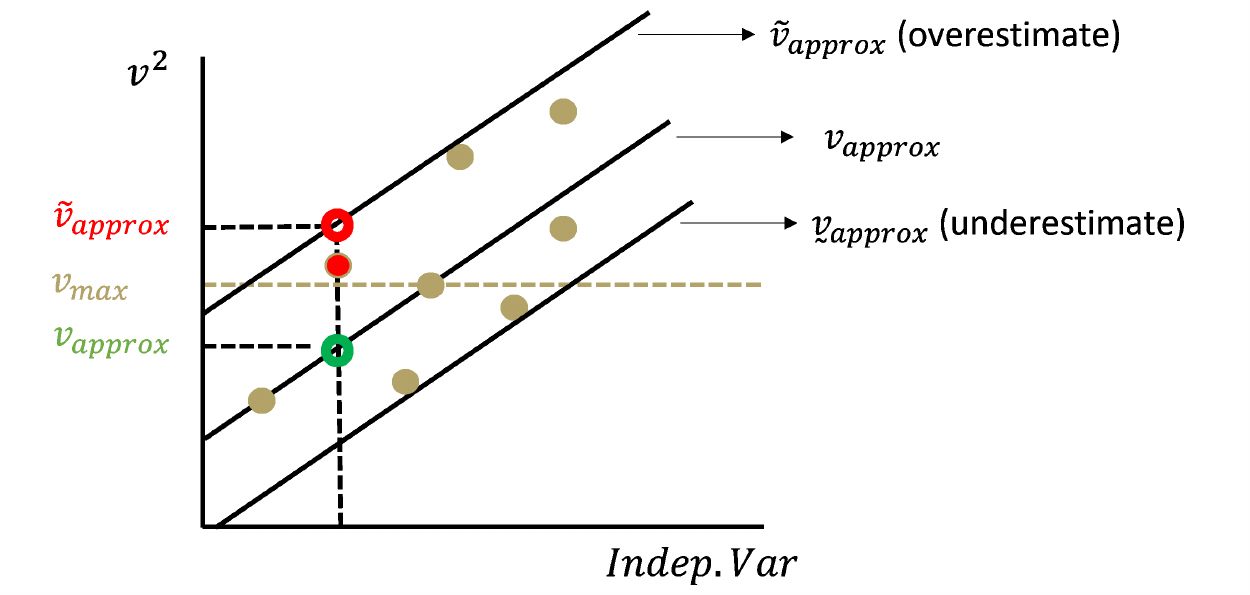}
    \caption{Conceptual illustration of the conservativeness of a pair of CLAs (over and underestimates). Using the pair of CLAs (top and bottom lines) would not lead to violations of lower and upper bounds imposed on $v$ as it would when using the non-conservative linear fit in the middle. In this example, unlike the overestimating CLA, the non-conservative approximation would erroneously predict satisfaction of the upper bound $v_{\sf max}$ for the red point.} 
    \label{fig: conservative_vs_nonconservative_approx}
    \vspace*{-1.2em}
\end{figure}

%
Thus, for the EEV-C problem, we model the distribution network through the sample-based linearization method first proposed by \cite{pap_2021} within the context of transmission networks. This methodology uses constrained regression problems to construct conservative linear approximations (CLAs) of target quantities. The approximations are said to be conservative since they are purposely constructed to either overestimate or underestimate the desired parameters, as exemplified in Fig.~\ref{fig: conservative_vs_nonconservative_approx}. The resulting affine equations can then be used to write upper and lower bounds in optimization problems.

To provide a time-series model of an entire distribution network, we propose constructing the following CLAs: 
\begin{subequations}
    \label{dist_CLA_approximated_parameters}
    \begin{align}
        &\widetilde{v}_i^t = \widetilde{g}_i(\vp_{\sf EV}^t) &\quad \forall i \in \setN, \forall t \in \setT  \label{v_upper}\\ 
        &\undertilde{v}_i^t = \undertilde{g}_i(\vp_{\sf EV}^t) &\quad \forall i \in \setN, \forall t \in \setT, \label{v_lower}
    \end{align}
\end{subequations}
where, unlike transmission network CLAs, the approximations $\widetilde{v}_i^t, \, \undertilde{v}_i^t \in \R$ represent overestimates and underestimates of the squared voltage magnitude of each \emph{single-phase node} $i \in \setN$ at time steps $t \in \setT$. These approximations are defined as affine functions of $\vp_{\sf EV}^t$, which denotes the vector of EV active power demands at time $t$ of buses $k \in \setK \subset \setN $ with designated EVs in them. We symbolize these functions as \emph{overestimating} and \emph{underestimating} CLAs $\widetilde{g}_i^t \: : \R^{|\setK|} \to \R$ and $\undertilde{g}_i^t \: : \R^{|\setK|} \to \R$, respectively.  
Note that each approximation is time specific, and that the background loads at each time step are treated as constants, \emph{not} input variables to these functions. Rather than voltage magnitudes, we construct approximations of their \emph{squares} to reduce linearization error,  as suggested in \cite{pap_2021}.

Consider the construction of one CLA $g_i^t$. Let $\va_{i,1}^{t} \in \R^{|\setK|}$, and $a_{i,0}^{t} \in \R$ be coefficients such that $g_i^{t}$ takes the form of an affine function. Let $\vp_{\sf EV}(m) \: \in \R^{|\setK|}$ denote vectors of $m = 1,\dots,M$ \emph{samples,} where each vector randomizes which EVs are charging in buses $k \in \setK \subset \setN$. Collect these sample vectors as columns in a sample matrix $\mP_{\sf EV} \in \R^{\abs{\setK} \times M}$. 
%
Furthermore, let  $\vv^{t}_i \: \in \R^{M}$ be a vector of computed target measurements of the parameter being approximated by $g_i^t$. Each entry in this vector is obtained by solving a power flow for each  sample vector  $\vp_{\sf EV}(m)$, and it takes the form
\begin{equation}
    \label{eq:vm2-samples}
    \vv^{t}_i \triangleq \begin{bmatrix}
        v_i^{t}(1) & \dots & v_i^{t}(m) & \dots v_i^{t}(M) 
    \end{bmatrix}^\top.
\end{equation}
Now we can construct an approximation of a specific squared nodal voltage $g_i^t$ at each time step $t\in\mathcal{T}$ via an $\ell_1$-norm approximation program of the form
\begin{subequations}
    \label{constrained_regression}
    \begin{align} 
        & \minimize_{a_{i,0}^t,\va_{i,1}^t} \norm{a_{i,0}^t \mathds{1} + \mP_{\sf EV}^\top \boldsymbol{a}_{i,1}^t - \vv_i^t}_1 \\
        &\quad \text{s.t} \quad a_{i,0}^t\mathds{1} + \mP_{\sf EV}^\top \boldsymbol{a}_{i,1}^t \leq \vv_i^t \quad \text{if underestimate} \label{if_under} \\
        &\quad \hphantom{\text{s.t}} \quad a_{i,0}^t\mathds{1} + \mP_{\sf EV}^\top \boldsymbol{a}_{i,1}^t \geq \vv_i^t \quad \text{if overestimate} \label{if_over},
    \end{align}
\end{subequations}
where $(\cdot)^\top$ is the matrix transpose and $\mathds{1}$ is an all-one vector.

If building an overestimate (resp. underestimate), constraints \eqref{if_over} (resp. \eqref{if_under}) are included to ensure that all predictions made by the resulting CLA are  above (resp. below) the sampled measurements. The $\ell_1$-norm objective can be motivated by empirical indications that the coefficients corresponding to linear approximations of distribution network models are often sparse. 
The program \eqref{constrained_regression} is a linear program with affine inequality constraints, which can handled with mature solvers.

Lastly, a time-varying distribution network model can be formulated within optimization problems using the previously computed CLAs via the following upper and lower bounds:
\begin{subequations}
    \label{dist-CLA-model}
    \begin{align}
        &\widetilde{g}_i^{t}(\vp_{\sf EV}^t) \leq \overline{v_i} &\quad \forall i \in \mathcal{N} , \forall t \in \mathcal{T} \label{dist-CLA-model_upper}\\ 
        &\undertilde{g}_i^{t}(\vp_{\sf EV}^t) \geq \underline{v_i} &\quad \forall i \in \mathcal{N} , \forall t \in \mathcal{T}, \label{dist-CLA-model_lower}
     \end{align}
\end{subequations}  
where $\overline{v_i}$ and $\underline{v_i}$ represent the maximum and minimum squared voltage magnitudes at node $i$, respectively. We next explain the need for \emph{conservative} approximations to obtain stronger constraints than simple linear regression, as in \eqref{dist-CLA-model}. We also explain how to solve the \mbox{EEV-C} problem while only explicitly enforcing a subset of these constraints.
\subsection{Trade-off between charging time and constraint violations}
\label{sec:tradeoff}
The EEV-C problem balances an inherent trade-off between the objective of minimizing total EV charging time and the desire to protect the power system from dangerous constraint violations. Given current distribution grid infrastructures, previous work on this topic in~\cite{nilsson_owen_aquino_coogan_molzahn-GreenEVT} has shown that network violations increase greatly with growing numbers of simultaneously charging EVs. Since the optimization algorithm will naturally try to simultaneously charge as many EVs as possible, grid constraints may render infeasible some of the most desirable solutions, thereby delaying 
the  evacuation process. 

For this reason and the urgency of an emergency evacuation, grid operators may want to allow the network to operate above its normal limits. To model this choice, we introduce the non-negative slack variables $\lambda^{t,+}_i$ and $\lambda^{t,-}_i$ for the upper and lower bounds in \eqref{dist-CLA-model}, respectively. Additionally, we add a new constraint \eqref{violations-upper-bound} to the CLA distribution network model that upper bounds the sum of all the slack variables. The  upper bound $\Lambda_{\sf max}$ allows the operator to pick a tolerable cumulative violation amount across the entire network. The CLA distribution model with allowable violations is then: 
\begin{subequations}
    \label{slacked-model}
    \begin{align}
        &\widetilde{g}_i^{t}(\vp_{\sf EV}^t) \leq v_{i, \mathsf{max}} + \lambda^{t,+}_i &\quad \forall i \in \mathcal{N} , \forall t \in \mathcal{T} \label{slacked-constraints_upper}\\ 
        &\undertilde{g}_i^{t}(\vp_{\sf EV}^t) \geq v_{i, \mathsf{min}} - \lambda^{t,-}_i &\quad \forall i \in \mathcal{N} , \forall t \in \mathcal{T} \label{slacked-constraints_lower}\\
        &\sum_{i \in \setN} \sum_{t \in \setT} (\lambda^{t,+}_i + \lambda^{t,-}_i) \leq \Lambda_{\mathsf{max}} \label{violations-upper-bound}\\
        &\lambda^{t,+}_i , \lambda^{t,-}_i \geq 0 &\quad \forall i \in \mathcal{N} , \forall t \in \mathcal{T}. \label{non-negative-slack}
     \end{align}
\end{subequations}  

\section{Iterative EEV-C Problem Using Conservative Linear Approximations (CLA-EEV-C) }\label{sec:CLA-EEV-C}

We formulate the EEV-C problem using the surrogate distribution network model described in Section \ref{sec:distribution_opt_models}. This grants the grid operator the flexibility to select a tolerable threshold of network violations throughout the system. Additionally, we introduce an iterative constraint generation algorithm designed to enhance tractability when solving the problem. 

\subsection{Formulation}
\label{sec:formulation}
We first state our modeling assumptions:
\begin{mdframed}[innerleftmargin=3pt, innerrightmargin=3pt]
    \noindent \underline{\textbf{Assumption 1:}} The background loads (without any EV charging) follow the patterns of a typical mid-summer day.
    
    \noindent \underline{\textbf{Assumption 2:}} EV loads operate at unity power factor, and their charge rate is 7.5 kW. It takes 32 time periods (8 hours) to charge an EV from 0\% to 100\% at this rate.
    
    \noindent \underline{\textbf{Assumption 3:}} When a TAZ is given the order to start charging, its EVs must charge until they reach full capacity. An EV may only charge after its TAZ was given the order to start charging.
    
    \noindent \underline{\textbf{Assumption 4:}} There are no more than 96 time periods (24 hours) available to charge all EVs. 

    \noindent \underline{\textbf{Assumption 5:}} The starting battery level of all EVs is treated as a known input parameter.

    \noindent \underline{\textbf{Assumption 6:}} The linearized distribution system model only considers nodal voltage violations.
\end{mdframed}
We note that including line current limits is the subject of ongoing investigation, to be covered in future work.

The objective of the EEV-C problem is to minimize the number of time periods, each of 15 minutes, that it takes to charge all EVs in a region while satisfying an upper bound on the magnitude of grid violations imposed by the grid operator. As a result, the algorithm pushes the charging of all vehicles as close as possible to their evacuation departure deadlines, thereby minimizing the amount of preparation time needed ahead of an incoming disaster.  The proposed model uses a well-studied framework for evacuation that partitions a city into Transportation Analysis Zones (TAZs). As evacuation instructions are given at the TAZ level (rather than feeder, street, or neighborhood levels), they are both realistic and computationally meaningful \cite{hasan_large-scale_part1_2021,hasan_large-scale_part2_2021}. Concretely, we denote the set of all TAZs as $\Xi$, where every $\xi \in \Xi$ is an individual TAZ with multiple EVs registered in it. Let $E_\xi = \{1,\dots,\abs{E_\xi}\}$ be  the set of electric vehicles registered within TAZ $\xi$, and $\mathcal{K}$ denote the set of buses $k \in \mathcal{K} \subset \setB$ with EVs in them.

Altogether, the EEV-C problem is presented in \eqref{ev-evac}. We denote $\Gamma_{\sf max} \in \R$ as the time step when the first TAZ starts charging. The objective \eqref{EEV-C-objective} maximizes $\Gamma_{\sf max}$ to make it as close as possible to the first moment of the evacuation. This is defined implicitly via \eqref{1st_eq} using $\tau^t \in \{0,1\}$, where $\tau^t =1$ if some TAZ has started charging prior to $t$ and 0 otherwise.

\begin{subequations}\label{ev-evac}
\noindent\hrulefill

\textbf{Emergency Electric Vehicle Charging (EEV-C) Problem.} 

\noindent\hrulefill
    \begin{align}
        &\: \: \: \: \maximize_{\boldsymbol{\tau},\boldsymbol{C}_{\xi,h},\boldsymbol{C}_\xi} \quad \Gamma_{\sf max}, \quad  \text{subject to:} \label{EEV-C-objective}\\
        &\: \: \: \: \Gamma_{\sf max} \leq t\tau^t + T(1-\tau^t) \quad \forall t \in \setT \label{1st_eq}\\
        &\: \: \: \:  \tau^t \geq \frac{\sum_{\xi \in \Xi} \sum_{t'=1}^{t} C_\xi^{t'}}{T |\Xi|} \quad \forall t \in \setT\label{2nd_eq}\\
        & \: \: \: \: L_{\xi,h}^t = L_{\xi,h}^{t-1} + \frac{1}{\beta}C_{\xi,h}^{t-1} \quad \forall t \in \setT, \forall \xi \in \Xi, \forall h \in E_\xi \label{3rd_eq}\\
        &\: \: \: \:  L_{\xi,h}^0 \leq L_{\xi,h}^{t} \leq 1 \quad \forall t \in \setT, \forall \xi \in \Xi, \forall h \in E_\xi \label{4th_eq}\\
        & \: \: \: \: C_\xi^t \geq C_\xi^{t-1} - \frac{\sum_{h \in E_{\xi}} L_{\xi,h}^t}{\abs{E_\xi}} \quad \forall t \in \setT, \forall \xi \in \Xi \label{5th_eq}\\
        & \: \: \: \: C_\xi^t \leq 2 - \frac{1+ \beta \sum_{h \in E_{\xi}} L_{\xi,h}^t}{\beta \abs{E_\xi}} \quad \forall t \in \setT, \forall \xi \in \Xi \label{6th_eq} \\
        &\: \: \: \:  C_\xi^{t} - L_{\xi,h}^t \leq C_{\xi,h}^t \leq C_\xi^t \quad \forall t \in \setT, \forall \xi \in \Xi, \forall h \in E_\xi \label{7th_eq}\\
        & \: \: \: \: \frac{\sum_{h \in E_{\xi}} L_{\xi,h}^{d_\xi}}{\abs{E_\xi}} = 1 \quad \forall \xi \in \Xi \label{8th_eq} \\
        & \: \: \: \: p_{EV_k}^t = \sum_{\xi \in \Xi} \sum_{h \in E_k} C_{\xi,h}^t R \quad \forall t \in \setT \label{power_transport} \\
        \nonumber & \: \: \: \:\eqref{slacked-constraints_upper}\text{--}\eqref{non-negative-slack} \\
        & \: \: \: \: \tau^t , C_\xi^t , C_{\xi,h}^t \in \{0,1\} \quad \forall t \in \setT, \forall \xi \in \Xi, \forall h \in E_\xi.
     \end{align}
     \noindent\hrulefill
\end{subequations} 

The main decision of interest in this program is \textit{when} each of the TAZs start charging. Once a TAZ $\xi \in \Xi$ starts charging, \textit{all} EVs $h \in E_{\xi}$ in that TAZ will charge until they reach full capacity. The charging progression of each EV is modeled by equations \eqref{3rd_eq}--\eqref{4th_eq}, where each EV $h \in E_{\xi}$ has an initial state of charge of $L_{\xi,h}^0$ and a battery level  $L_{\xi,h}^t$ at time $t \in \mathcal{T}$. Additionally, $\beta$ denotes the number of time steps it takes to fully charge an EV from 0\% at rate $R$, and each $\xi \in \Xi$ is assigned a binary variable $C_\xi^t \in \{0,1\}$, which takes the value of 1 if TAZ $\xi$ is charging at time $t$, and 0 otherwise. After all EVs in a TAZ are at full capacity, the entire TAZ will stop charging. This behavior is modeled by constraints \eqref{5th_eq}--\eqref{6th_eq}, where the binary variables $C_{\xi,h}^t \in \{0,1\}$ take the value of 1 if vehicle $h$ in TAZ $\xi$ is charging at time $t$, and 0 otherwise. Furthermore, constraint \eqref{7th_eq} ensures that an EV can only charge if its TAZ is instructed to charge.

Constraint \eqref{8th_eq} ensures that each EV in a TAZ $\xi$ is fully charged by its TAZ departure time $d_\xi$. Constraint \eqref{power_transport} links the EV demand across time to their charging schedules. Constraints \eqref{slacked-constraints_upper}--\eqref{non-negative-slack} represent the \emph{surrogate} network model, and its linearized network constraints described in Section~\ref{sec:distribution_opt_models}. Since the surrogate network constraints (as well as all other constraints related to the charging of EVs) are linear, we have a mixed-integer linear program (MILP).

\subsection{Solution algorithm}
\label{sec:algorithm}
\begin{figure}[t!]
    \centering
    \includegraphics[width=0.801628\linewidth]{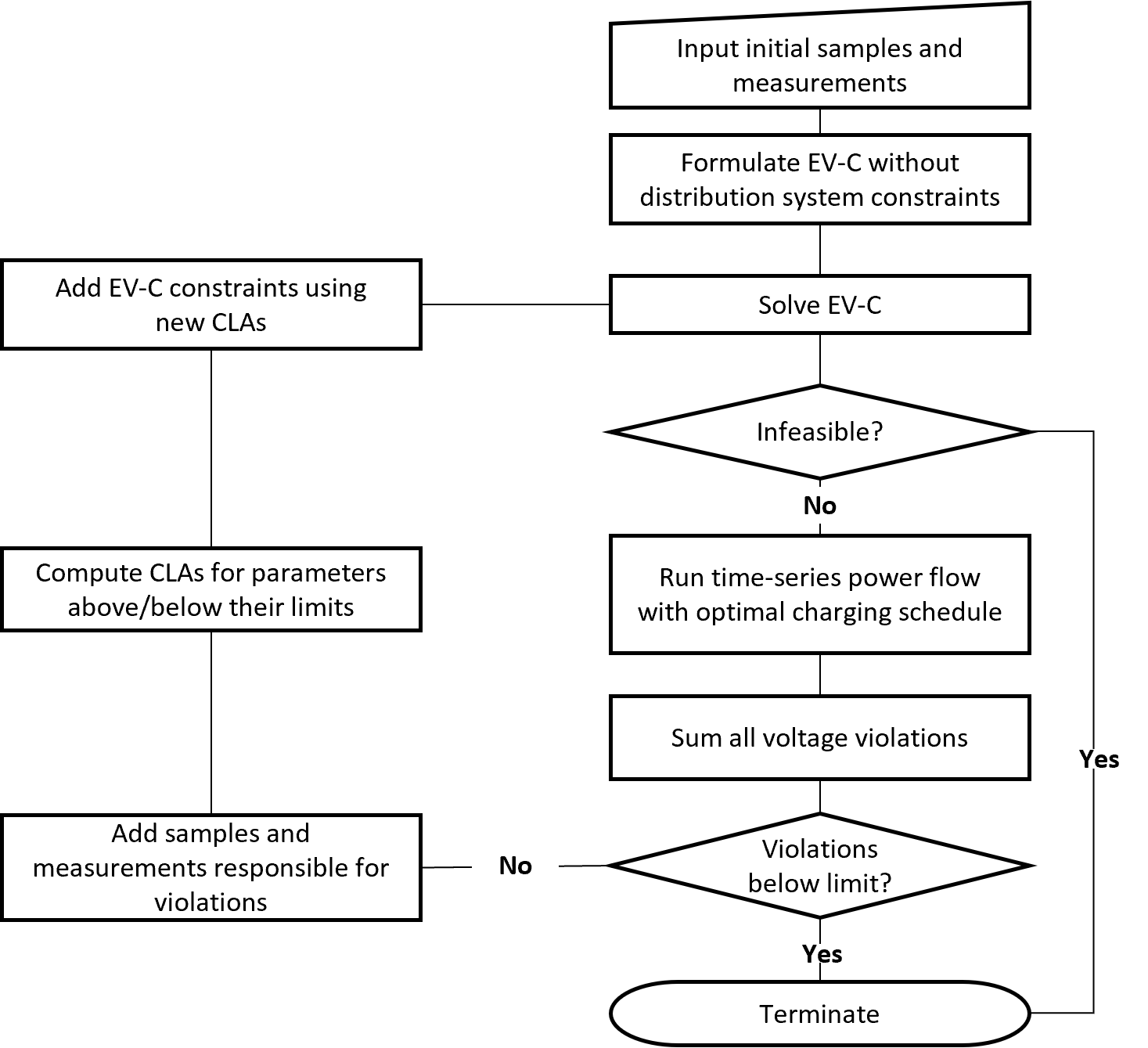}
    \caption{Constraint generation flowchart for the EEV-C problem.} 
    \label{fig: flowchart}
    \vspace*{-0.8em}
\end{figure}

Applying the constraints \eqref{slacked-constraints_upper}--\eqref{slacked-constraints_lower} to the EEV-C problem \eqref{ev-evac} provides the advantage of simplicity relative to the AC OPF constraints. Nonetheless, the constraints can still render an intractable model, particularly for large  networks. Even when taking advantage of parallel computing, generating all of these CLAs would require significant up front computational effort. For these reasons, we propose solving \eqref{ev-evac} using the iterative constraint generation scheme outlined in Fig. \ref{fig: flowchart}.
%

The algorithm begins by processing an initial matrix $\mP_{\sf EV}$ and vectors $\vv^{t}_i$ for all $(i,t) \in \setN \times \setT$. First, the algorithm solves the EEV-C problem \eqref{ev-evac} without voltage constraints \eqref{slacked-constraints_upper}-\eqref{slacked-constraints_lower}, generating a ``naive" charging schedule that disregards potential power system violations. Subsequent time-series power flow analysis on this schedule identifies any actual voltage violations. If these exceed the grid operator's allowable limit $\Lambda_{\sf max}$, the algorithm iteratively adds constraints until the violation limit is met or infeasibility is proven.

New constraints are formed by constructing CLAs through constrained regression problems like \eqref{constrained_regression}. These constraints are based on parameters that violated limits in the previous iteration's charging schedule, as assessed by a power flow solver. With each iteration, measurements of these violated parameters and the corresponding EV active power demands $\vp_{\sf EV}^t$ are added to the sets of measurements and samples.

Using conservative approximations significantly speeds up the convergence of the iterative approach. Since the conservative constraints are stronger than those created from simple least-squares regression, they exert greater influence on forcing new solutions to the EEV-C problem at each iteration. This, in turn, accelerates the algorithm's convergence, leading to fewer iterations and reduced computation time.  

\section{Numerical Experiments and Discussion}
\label{sec:numerical-experiments}
We demonstrate the proposed methods on a synthetic model of the distribution network of the city of Greensboro, NC. To achieve this, we make use of the  Greensboro Electric Vehicle Testbed (\texttt{GreenEVT})\footnote{\texttt{GreenEVT} is available at \url{https://github.com/GreenEVT/GreenEVT}.}, which couples each bus in the power network to its corresponding Transportation Analysis Zone (TAZ) in the transportation system. This testbed is built on top of NREL's \texttt{SMART-DS} (Synthetic Models for Advanced, Realistic Testing: Distribution systems and Scenarios) dataset~\cite{palmintier_experiences_2021}, which provides realistic-but-not-real distribution network datasets for three regions (industrial, rural, and urban-suburban) for the city of Greensboro. The \texttt{GreenEVT} testbed provides four EV penetration scenarios: low, medium, high, and extreme. The time-series load data models the entire time needed to charge all EVs before an evacuation. 

For this paper's test case, we will simulate the high EV penetration scenario for substation 19 in the urban-suburban region, which contains 3 feeders and a total of 4815 single phase nodes. All optimization problems, including the \mbox{EEV-C} problem and the constrained regression problems used to compute CLAs, were solved using Gurobi v10.0.1. Additionally,  all the power flow solutions were computed using the power flow simulator \texttt{OpenDSS}~\cite{dugan_opendss_2011} and the Python interface \texttt{yadi}~\cite{yadi_github}. The computations were carried out on Georgia Tech's PACE cluster using a computing node equipped with a quad-core 2.7 GHz processor and 64 GB of RAM.

\begin{figure}
    \centering
    \includegraphics[width=0.8\linewidth]{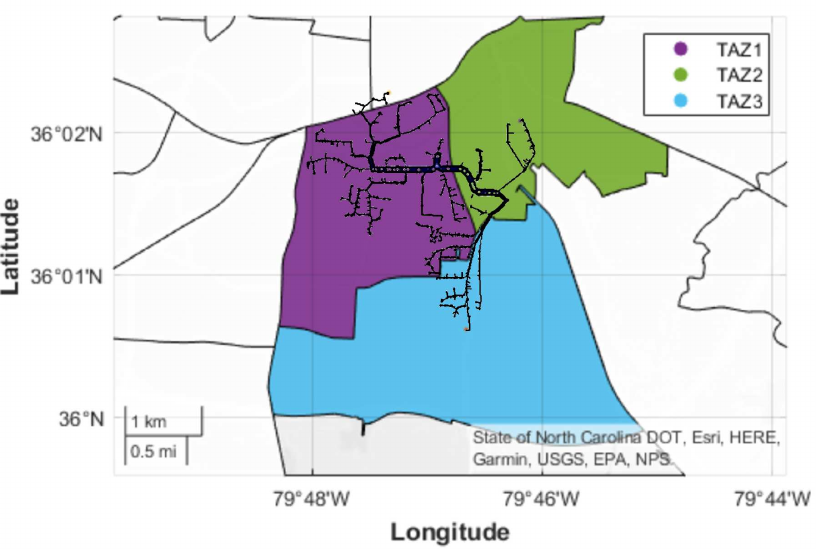}
    \caption{Geographic map of the distribution network downstream of substation 19 in Greensboro's urban-suburban network. TAZ1, TAZ2, and TAZ3, which surround this network, are shown in blue, red, and yellow, respectively.} 
    \label{fig: geoplot}
    \vspace*{-1.2em}
\end{figure}

\subsection{Total charging time vs. Allowable violations}
\label{sec:allowable_violations_vs_charging_time}
\begin{figure}[t]
    \centering
    \includegraphics[width=0.84\linewidth]{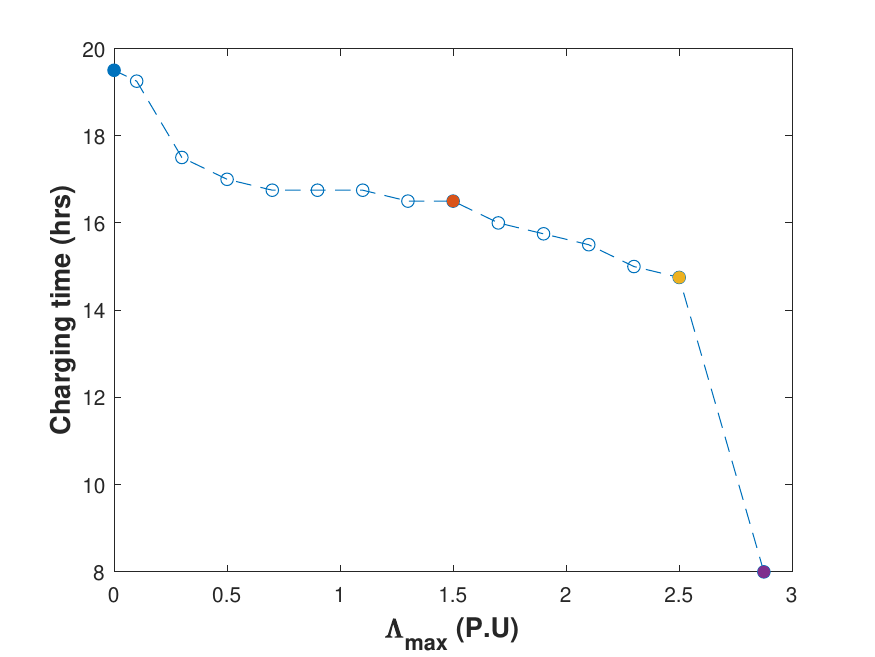}
    \includegraphics[width=0.84\linewidth]{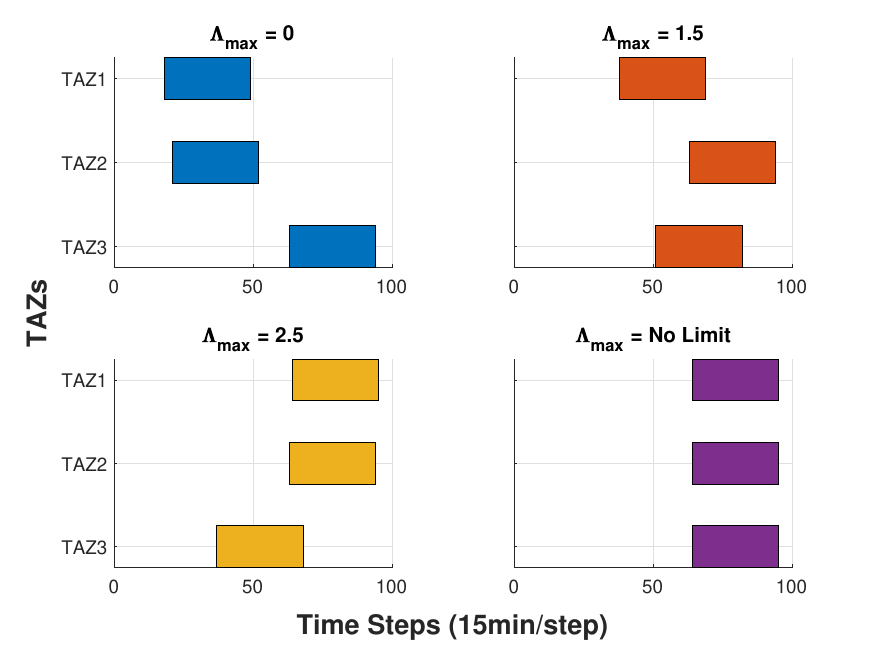}
    \caption{The top graph shows the optimal charging time at varying allowable violation thresholds. The bottom four graphs show the charging schedules 
    that achieve the optimal charging times at specific points in the top graph.} 
    \label{fig: allowable_violations_vs_charging_time}
\end{figure}

Solving the EEV-C problem across multiple values for $\Lambda_{\sf max}$ reveals the trade-off between total charging time and network violations, shown in the top of Fig.~\ref{fig: allowable_violations_vs_charging_time}. As expected, an increase in permissible violations corresponds with a reduction in the time required to charge the EVs in all three TAZs in the studied test case. Note that significant reductions in charging time occur primarily at the lower and upper ends of this curve, featuring a 9.09\% decrease when increasing $\Lambda_{\sf max}$ from 0 to 0.4, and a substantial 45.76\% decrease when increasing $\Lambda_{\sf max}$ from 2.5 to 2.87. Note that raising $\Lambda_{\sf max}$ to 2.87 mirrors the absence of constraints on the distribution network, resulting in the ``naive'' charging schedule. 

The bottom of Fig.~\ref{fig: allowable_violations_vs_charging_time} uses Gantt charts to visualize the optimal charging schedules obtained under different values of $\Lambda_{\sf max}$. In these charts, each horizontal bar depicts the time intervals for charging each TAZ whose vehicles are all scheduled to depart at $t = 96$.  
Interestingly, when $\Lambda_{\sf max} = 0$, a 2.75 hour gap appears between time steps 52 and 63, during which no vehicles are actively charging. This behavior suggests that background loads at these time steps are causing considerable delays in the charging of TAZ1 and TAZ2. This indicates that demand response strategies to prioritize critical EV charging loads prior to evacuations may be advantageous.

\subsection{Analysis of the iterative algorithm}
\label{sec: analysis_of_iterative_algorithm}
Each data point in the upper graph of Fig.~\ref{fig: allowable_violations_vs_charging_time} is obtained using the iterative algorithm detailed in Section~\ref{sec:algorithm}. For this specific test case, all data points converged within an average time of 4.01 hours, and in between 1 and 4 iterations Additionally, these points converged with an average relative error of 0.68\% between predicted and actual cumulative violation magnitudes in the last iteration. This result points to the strong predictive performance of the proposed CLA distribution network model, particularly in the last iterations of the algorithm, where a larger set of CLAs constraints are used to describe the network. 

\begin{figure}
    \centering
    \includegraphics[width=0.84\linewidth,keepaspectratio]{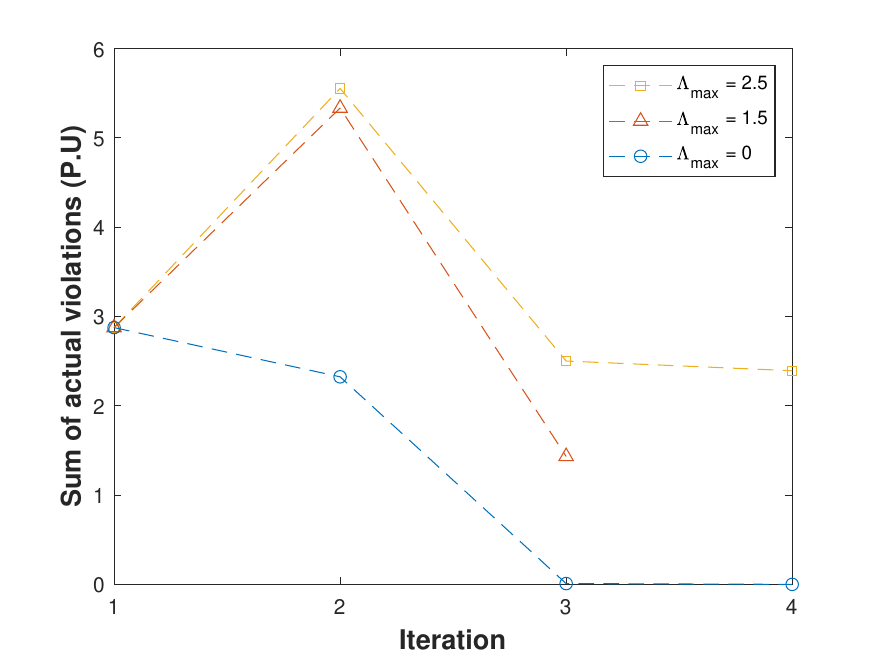}
    \includegraphics[width=0.84\linewidth,keepaspectratio]{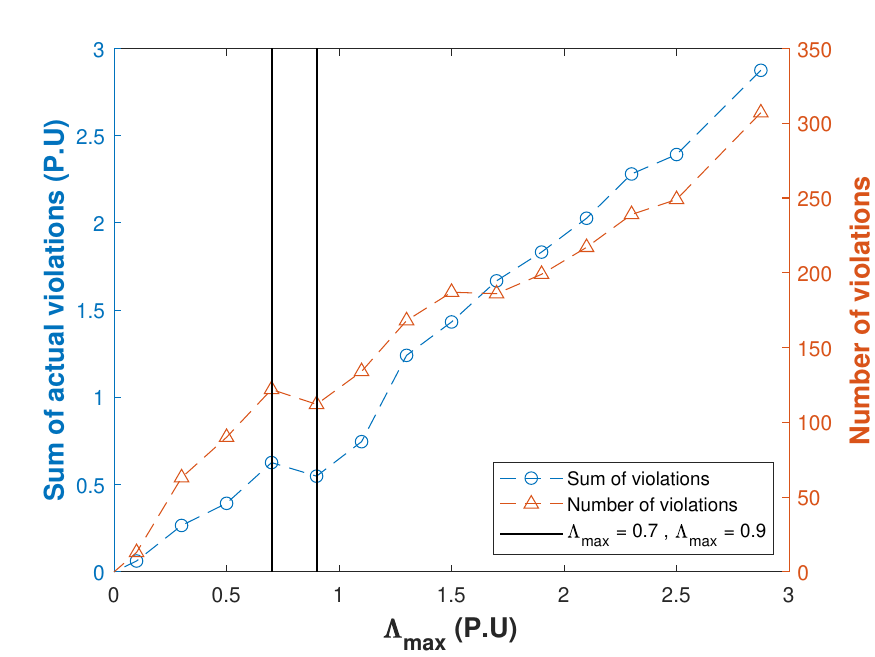}
    \caption{The top graph shows the trajectory of the sum of actual violations during the iterative algorithm for three different values of $\Lambda_{\sf max}$. The bottom graph shows the sum of violation magnitudes and the number of violations vs. $\Lambda_{\sf max}$ using the left and right y-axes, respectively. The two vertical lines at $\Lambda_{\sf max} = 0.7$ and  $\Lambda_{\sf max} = 0.9$ show the curves' non-monotonic nature.} 
    \label{fig: iterative_algorithm_analysis}
    \vspace*{-1.2em}
\end{figure}

The iterative process of obtaining the EEV-C solution when $\Lambda_{\sf max} = 0$, $\Lambda_{\sf max} = 1.5$, and $\Lambda_{\sf max} = 2.5$ is shown in the top graph of Fig.~\ref{fig: iterative_algorithm_analysis}. At each iteration, the algorithm first runs the EEV-C optimization problem and then simulates the resulting charging schedule in OpenDSS, so  each run converges when the sum of \emph{actual} violations reaches the respective $\Lambda_{\sf max}$ threshold. Since the algorithm always starts by solving the problem without any grid constraints, all three curves share a common starting point at the first iteration. Subsequent iterations reveal how the algorithm's approach does not always result in an immediate reduction of actual violations, but successfully converges after a few iterations.   

The proposed solution algorithm and EEV-C formulation produces monotonically decreasing charging times over $\Lambda_{\sf max}$, as anticipated. However, note that it does not produce a monotonically increasing sum of violation magnitudes nor violation counts, as one would expect. This behavior is prominently depicted in the bottom graph of Fig.~\ref{fig: iterative_algorithm_analysis}, where an increase in $\Lambda_{\sf max}$ from 0.7 to 0.9 results in a \emph{decrease} in both the total violation count and the cumulative violation magnitudes. This observation indicates that once the algorithm identifies an optimal charging schedule for a particular $\Lambda_{\sf max}$, it does not continue to explore solutions that maintain the same charging time while further reducing network violations.

\section{Conclusion}
\label{sec:conclusion}
In this work, we proposed an emergency EV charging scheduling formulation that leverages an embedded distribution network model featuring conservative power flow linearizations. The formulation provides flexibility by enabling the user to define an acceptable constraint violation threshold to  strike a balance between minimizing total EV charging time and safeguarding the power system from potential hazards. To tackle the computational complexity of the problem, we also presented an algorithm that explicitly enforces only a subset of the linearized power grid constraints at each iteration, thereby contributing to a more manageable solution approach.

The proposed algorithm was tested on a small region of a larger distribution network model of Greensboro, NC, under a high EV penetration scenario. The results demonstrate the effectiveness and accuracy of the linearized distribution model used in the formulation, and provide an illustration of the marginal utility of allowing different degrees of cumulative voltage violations in order to reduce the EV charging timeline. However, the experiments reveal that while the solution algorithm successfully minimizes total EV charging time within the defined violation threshold, it does not always minimize violations while maintaining the optimal charging time. These results motivate directions for future research:
\begin{enumerate}
    \item Extending the proposed formulation to co-minimize charging times and network violations.
    \item Exploring decomposition techniques to scale the formulation for the entire city of Greensboro, including its 21 substations and over 60 TAZs.
    \item Incorporating line and transformer current CLAs into the distribution network model and assessing the tractability of the resulting formulations.
    \item Integrating this algorithm into a citywide evacuation plan via a departure-scheduling-and-routing problem.
\end{enumerate}


\section*{Acknowledgement}
This work was supported in part by the Strategic Energy Institute at Georgia Tech, the NSF Graduate Research Fellowship Program under Grant No. DGE-1650044, and NSF award \#2145564

\bibliographystyle{ieeetr}
\bibliography{references}

\end{document}